\pdfoutput=1
\documentclass[preprint]{article}

\usepackage{amsmath,amsthm,enumitem,multicol,amssymb,mathrsfs,verbatim,graphicx}
\usepackage{geometry,xcolor}
\usepackage[makeroom]{cancel}
\usepackage[normalem]{ulem}
\makeatletter
\def\uwave{\bgroup \markoverwith{\lower3.5\p@\hbox{\sixly \textcolor{red}{\char58}}}\ULon}
\font\sixly=lasy6
\makeatother

\DeclareMathSymbol{\R}{\mathbin}{AMSb}{"52}
\newcommand{\A}{\mathcal{A}}
\newcommand{\B}{\mathscr{B}}
\newcommand{\C}{\mathbb{C}}

\newcommand{\E}{\mathbb{E}}

\renewcommand{\R}{\mathbb{R}}

\newcommand{\cross}{\text{Cr}}

\newcommand{\Cr}{\text{Cr}}
\newcommand{\CNB}{\text{CNB}}
\newcommand{\Nest}{\text{Nest}}

\DeclareMathSymbol{\N}{\mathbin}{AMSb}{"4E}
\DeclareMathSymbol{\Z}{\mathbin}{AMSb}{"5A}

\renewcommand{\P}{\mathcal{P}}

\newtheorem*{assump}{Main Hypotheses}
\newtheorem{remark}{Remark}

\newtheorem{theorem}{Theorem}

\begin{document}

 \title{Type B Gaussian Statistics as Noncommutative Central Limits}

\author{Natasha Blitvi\'c\thanks{Corresponding author.}\\Department of Mathematics and Statistics, Lancaster University\\ Lancaster, LA1 4YW, United Kingdom\\ \texttt{natasha.blitvic@lancaster.ac.uk}\\ \and  Wiktor Ejsmont\thanks{Partially supported by the National Science Centre grants no.\ 2014/15/B/ST1/00064 and no.\ 2018/29/B/HS4/01420.}\\Mathematical Institute University of Wroclaw \\
pl. Grunwaldzki 2/4, 50-384 Wroclaw, Poland\\
\texttt{wiktor.ejsmont@math.uni.wroc.pl}}

\maketitle

\begin{abstract}
We show that the noncommutative Central Limit Theorem of Speicher can be adapted to produce the Gaussian statistics associated to Coxeter groups of type B, in the sense of Bo\.zejko, Ejsmont, and Hasebe. Viewed through the lens of central limits, the passage from $q$-Gaussian statistics, associated with symmetric groups, to the Gaussian statistics associated with Coxeter groups of type B is precisely the passage from a sequence of independent elements that pairwise commute or anticommute, to a coupled pair of such sequences. The results pave the way for the transfer of known results from the bosonic/fermionic settings to such broader contexts.
\end{abstract}

\section{Introduction}

\emph{Noncommutative probability} is broadly concerned with the (noncommutative) distributions of objects arising from algebraic or operator algebraic contexts. Outside of the rich setting of quantum probability (see e.g. \cite{meyer2013quantum}) and the celebrated free probability of Voiculescu \cite{Voiculescu1992}, there exists a number of noncommutative probabilistic frameworks, such as the `$q$-deformed' probability \cite{Bozejko1991,Zagier1992,Bozejko1997}, that mirror to varying extents the central ideas of classical probability theory. Despite their seemingly specialist nature, these exhibit far-reaching connections to other areas of pure and applied mathematics, including combinatorics,
$q$-series, and physics (see e.g. [Bli12] and the references therein) and are therefore of broader mathematical interest.
A recent addition to this body of work is the generalized Gaussian process arising from Coxeter groups of type B, introduced by Bo\.zejko, Ejsmont, and Hasebe in \cite{Bozejko15}.

While quantum probability is grounded in physical reality through its ties to the bosonic/fermionic frameworks, the free probability naturally captures the scaling limits of large random matrices (see e.g.\ \cite{Biane2003}), and the $q$-deformed probability can be traced back to questions in quantum optics (see the review article \cite{Dodonov2002}), the idea of a generalized Gaussian process arising from Coxeter groups of type B may, at first sight, appear significantly more abstract and perhaps also farther removed from classical probabilistic intuition.

On the contrary, we show that this generalized Gaussian process naturally occurs in systems of `mixed spins'. Namely, through a new construction that draws on Speicher's noncommutative take on the classical Central Limit Theorem \cite{Speicher1992}, the noncommutative Gaussian distribution arising from Coxeter groups of type B is the central limit for ensembles of `independent' elements that pairwise commute or anticommute. In this sense, similarly to \cite{Bozejko1991,Speicher1992} and \cite{Blitvic2012,BlitvicCLT},  these generalized Gaussian statistics are a mixture of bosonic and fermionic statistics. Furthermore, under the lens of the Central Limit Theorem, the passage from $q$-Gaussian statistics, associated with symmetric groups (in the sense made clear in the following section), to the Gaussian statistics associated with Coxeter groups of type B is precisely the passage from a sequence of independent elements that pairwise commute or anticommute, to a coupled pair of such sequences. The construction is surprisingly elegant and paves the way to transferring known results from the bosonic/fermionic settings to the case at hand.

Prior to formulating our results, we begin by outlining the construction of the Gaussian statistics of `type B' and of a prototypical noncommutative Central Limit Theorem.

\subsection{Noncommutative Probability and the Type B Gaussian Statistics}

The results in this paper take place in a \emph{noncommutative probability space} $(\A,\varphi)$. The latter is formed by a unital $\ast$-algebra $\A$, whose elements are to be interpreted as `noncommutative random variables', and a state $\varphi$ on $\mathcal{A}$ (that is a linear functional on $\mathcal{A}$ satisfying $\varphi( a^* a) \geq 0$ for all $a \in \A$ and $\varphi(1) = 1$), playing the role of classical expectation. The \emph{distribution} of a noncommutative random variable $a\in\A$ is encoded by its \emph{mixed moments}, i.e. expressions of the form  $\varphi(a^{\epsilon(1)}\ldots a^{\epsilon(k)})$ for all $k\in\mathbb N$ and $\epsilon(1),\ldots,\epsilon(k)\in\{1,\ast\}$. The distribution of a self-adjoint element $a=a^\ast$ is encoded by its \emph{moment sequence} $(\varphi(a^k))_{k\in\mathbb N}$. 

The type B Gaussian elements of \cite{Bozejko15} `live' in the algebra of bounded linear operators on a deformed Fock space. We briefly outline their construction and its context. Recall that the \emph{Coxeter group of type B} (also known as \emph{hyperoctahedral group}) of degree $n$, denoted by $\Sigma(n)$, is generated by the elements $\pi_0,\pi_1,\ldots,\pi_{n-1}$ subject to the defining relations
\begin{eqnarray*}&\pi_i^2=1\text{ for } 0\leq i\leq n-1,\quad\quad & \pi_i\pi_j=\pi_j\pi_i\text{ for } 0\leq i,j\leq n-1\text{ s.t. } |i-j|>1,\\ &\pi_0\pi_1\pi_0\pi_1=\pi_1\pi_0\pi_1\pi_0,\quad\quad & \pi_i\pi_{i+1}\pi_i=\pi_{i+1}\pi_i\pi_{i+1}\text{ for } 1\leq i\leq n-2.
\end{eqnarray*}
Equivalently, $\Sigma(n)$ is the semidirect product $S(n) \ltimes \mathbb{Z}^n_2$ for the obvious action of the symmetric group $S(n)$ on the group $\mathbb{Z}^n_2$, or the wreath product $\mathbb{Z}_2 \wr S(n)$ for the natural action of $S(n)$ on the set $[n] := \{ 1, \ldots , n \}$.
More concretely, setting
$[n]_{\pm} := \{ \pm 1, \cdots , \pm n \}$,
$\Sigma(n)$ is the subgroup of the permutation group of $[n]_{\pm}$ consisting of
those elements that commute with the permutation $\tau$ given by inversion $k \mapsto -k$ ($k \in [n]_{\pm}$).
As generators of $\Sigma(n)$, one may now take
$\pi_0$ to be the transposition $(-1,1)$ and, for $1 \leq i \leq n-1$, $\pi_i$ to be the product of transpositions $(i, i+1)(-i, -i-1)$.
For each element $\sigma$ of $\Sigma(n)$ there is
a unique non-negative integer $k = k(\sigma)$ such that $\sigma = \pi_{i_1} \cdots \pi_{i_k}$ where
$0 \leq i_1, \cdots , i_k \leq n-1$ and $\sigma$ cannot be expressed as such a composition of less than $k$ generators. Moreover
the following quantities do not depend on the choice of such minimal representation (see \cite{Bozejko1994}, Theorem 2.1):
$$\ell_0(\sigma) := \# \big\{ p\in \{ 1, \ldots , k(\sigma) \}\mid i_p = 0 \big\}
\text{ and }
\ell(\sigma) := \# \big\{ p\in \{ 1, \ldots , k(\sigma) \}\mid i_p \neq 0 \big\}.$$
To view
the symmetric group $S(n)$ as a subgroup of $\Sigma(n)$,
identify $\lambda \in S(n)$ with
the element $\widehat{\lambda} \in \Sigma(n)$ defined by
$$\widehat{\lambda}(i) :=\begin{cases}\lambda(i) & i \in \{ 1, \cdots ,n \} \\
        -\lambda(-i)  & i \in \{ -1, \cdots , -n \}
         \end{cases}.$$
Thus, for $1 \leq i \leq n-1$,
$\pi_i = \widehat{\lambda}$ where $\lambda$ is the transposition $(i, i+1)$.

Next, let $H$ be a complex separable Hilbert space and consider an involutive unitary operator $\Pi_0$ on $H$ so that $\Pi_0$ is self-adjoint and $(\Pi_0)^2 = I_H$. Then, for each $n \in \mathbb{N}$, the group $\Sigma(n)$ acts unitarily on the Hilbert space $H^{\otimes n}$ as follows:
$$U_{\pi_0} := \Pi_0 \otimes I_H^{\otimes (n-1)}
\ \text{ and, for }
\lambda \in S(n), \ U_{\widehat{\lambda}} ( x_1 \otimes \cdots \otimes x_n ) := x_{\lambda(1)} \otimes \dots \otimes x_{\lambda(n)}.$$
Denote by $\mathcal F_{\text{fin}}(H)$ the  \emph{algebraic full Fock space}, namely, the algebraic direct sum $\bigoplus_{n\geq 0} H^{\otimes n}$ with 
$H^{\otimes 0} := \mathbb{C}$.
For $\alpha, q \in [-1,1]$, define the \emph{symmetrization operator} on $\mathcal F_{\text{fin}}(H)$ by
$P_{\alpha, q} := \bigoplus_{n \geq 0} P_{\alpha, q}^{(n)}$, where $P^{(0)}_{\alpha,q} := I_{\mathbb{C}}$ and, for $n \in \mathbb{N}$,
$$P^{(n)}_{\alpha,q} :=
\sum_{\sigma \in \Sigma(n)} \alpha^{\ell_0(\sigma)} q^{\ell(\sigma)} U_\sigma.$$
Then $P_{\alpha, q}$ is positive semi-definite, i.e.
$$\langle \xi, P_{\alpha, q} \xi \rangle = \sum_{n \geq 0} \langle \xi_n, P^{(n)}_{\alpha,q} \xi_n \rangle \geq 0 \quad \text{for } \xi = ( \xi_n) \in \mathcal{F}_{\text{fin}}(H),$$
with strict inequality for $\xi \neq 0$ if $|\alpha|, |q| < 1$ (\cite{Bozejko1994}, Theorem 2.1).
Therefore,
\[
\langle \zeta, \eta \rangle_{\alpha, q} :=
\big\langle \zeta, P_{\alpha, q} \eta \big\rangle.
\]
is an inner product on $\mathcal F_{\text{fin}}(H)$ for $\alpha,q\in (-1,1)$. The completion of the pre-Hilbert space $(\mathcal F_{\text{fin}}(H),\langle\,,\,\rangle_{\alpha,q})$ is denoted $\mathcal F_{\alpha,q}(H)$ and termed the \emph{$(\alpha,q)$-Fock space} or \emph{Fock space of type B}.

The Fock space of type B is a generalization of the $q$-Fock space of Bo\.zejko and Speicher \cite{Bozejko1991}, the latter identified with the case $\alpha=0$, for which many interesting probabilistic results are known (see e.g. \cite{Anshelevich2001,Biane1997b,Donati-Martin2003,Kemp2005,Anshelevich2010,Deya2013,Guionnet2014}). In the setting of Fock spaces, the relevant probabilistic aspects manifest through a family of operators that play the role of `noncommutative Gaussian' random variables. Specifically, for the Fock space of type B, the operators of interest are the following.

Fix $\alpha$ and $q$ in the interval $(-1,1)$. The $(\alpha, q)$-\emph{creation operator with test vector} $x \in H$ is defined first on the dense subspace $\mathcal{F}_{\text{fin}}(H)$ by 
$$b_{\alpha, q}^* (x) \xi :=
( \xi_n \otimes x )_n = ( 0, \xi_0 x, \xi_1 \otimes x, \cdots )\quad \text{ for } \xi = ( \xi_n ) \in \mathcal{F}_{\text{fin}}(H),$$
and then extended to an operator on $\mathcal F_{\alpha,q}(H)$ by continuity. (Its boundedness is proved in \cite{Bozejko15}, Theorem 2.9.) Its adjoint is 
called the $(\alpha, q)$-\emph{annihilation operator with test vector} $x$ and denoted $b_{\alpha, q} (x)$. The $(\alpha, q)$-\emph{Gaussian} operators are then defined by
$G_{\alpha, q}(x) := b_{\alpha, q}(x) + b_{\alpha, q}^*(x)$,
generalizing the $q$-Gaussian operators of \cite{Bozejko1991}. Use of the right creation operator (as opposed to the more usual left creation) gives more transparent compatibility with the symmetrizing operators $P^{(n)}_{\alpha, q}$ and the subgroup embeddings $\Sigma(n-1) \to \Sigma(n)$.

The distribution $\mu_{\alpha, q, x}$ of the $(\alpha, q)$-Gaussian $G_{\alpha, q}(x)$ in the vacuum state $T \mapsto \langle \Omega, T \Omega \rangle$, where $\Omega := ( 1, 0, 0, \cdots ) \in \mathcal{F}_{\alpha,q}(H)$, is the orthogonalizing probability measure of the $q$-Meixner--Pollaczek polynomials $MP_{\alpha \langle x, \Pi_0 x \rangle, q}$ (\cite{Bozejko15}, Theorem 3.3). As previously noted, the case $(\alpha=0,q)$ reduces to the $q$-Gaussian measure, while $(\alpha,q=0)$ recovers the symmetric free Meixner laws \cite{Saitoh2000,Anshelevich2001}.

The above operators satisfy the $(\alpha, q)$-commutation relations (\cite{Bozejko15}, Proposition 2.6):
\begin{equation}
b_{\alpha,q}(x)b^\ast_{\alpha,q}(y)-qb^\ast_{\alpha,q}(y)b_{\alpha,q}(x)=\langle x,y\rangle I+\alpha\langle x,\Pi_0 y\rangle q^{2N}\quad\quad (x,y\in H),\label{eq-comm}
\end{equation}
in which $I$ denotes the identity operator on $\mathcal{F}_{\alpha, q}(H)$ and $q^{2N}$ is the contraction operator on $\mathcal{F}_{\alpha, q}(H)$ given by $\xi = (\xi_n) \mapsto ( q^{2n} \xi_n )$.

Corresponding Gaussian processes arise by fixing a conjugation on $H$, so that $H$ is the complexification of a real Hilbert space $H_{\mathbb{R}}$. One way of doing this is by fixing an orthonormal basis $( e_n )$ for $H$ and letting $H_{\mathbb{R}}$ be the closed real-linear span of the basis; the conjugation is then given by $\sum z_n e_n \mapsto \sum \overline{z_n} e_n$.

We now fix a conjugation $x \mapsto \overline{x}$ on $H$ and
set $H_{\mathbb{R}} := \{ x \in H: \overline{x} = x \}$.
The corresponding type B Gaussian process,
$\big( G_{\alpha, q}(x) = b_{\alpha, q}(x) + b^*_{\alpha, q}(x) \big)_{x \in H_{\mathbb{R}}}$,
has moments in the vacuum state expressible via a Wick-type formula
\begin{align}
&\langle \Omega,G_{\alpha,q}(x_{2n-1})\ldots G_{\alpha,q}(x_1)\Omega\rangle_{\alpha,q}=0\\
&\langle \Omega,G_{\alpha,q}(x_{2n})\ldots G_{\alpha,q}(x_1)\Omega\rangle_{\alpha,q}=\notag\\& \sum_{(\pi,f)\in \P_2^B(2n)}\alpha^{\text{NB}(\pi,f)}\,q^{\cross(\pi)+2\text{CNB}(\pi,f)}\prod_{\substack{(i,j)\in\pi\\f(i,j)=1}}\langle x_i,x_j\rangle  \prod_{\substack{(i,j)\in\pi\\f(i,j)=-1}}\langle x_i,\Pi_0\, x_j\rangle
 \label{eq-moments}
\end{align}
where $x_1,\ldots,x_n\in H_{\R}$, $\P_2^B(2n)$ denotes the set of pair partitions of type B on $\{1,\ldots,2n\}$, $NB(\pi,f)$ is the number of negative blocks of $(\pi,f)$, and Cr and CNB are the crossing and asymmetric nesting statistics defined further in this paper (see Section~\ref{sec-comb}).

\subsection{Noncommutative Central Limit Theorems}
Our starting observation is that the form of \eqref{eq-moments} is evocative: by appearing as a product of covariances multiplied by a combinatorial statistic, the moments of the Gaussian operator of type B hint at the existence of a (noncommutative) Central Limit Theorem, whose combinatorial proof brings to the fore such pairwise structure.

Indeed, in \cite{Speicher1992}, Speicher showed that the $q$-Gaussian statistics ($\alpha=0$ case) arise from a noncommutative Central Limit Theorem. Speicher's central argument considers a sequence of elements with `mixed spins'. In its simplest form, the theorem concerns a sequence of self-adjoint elements $(a_i)_{i\in\mathbb N}$ of a noncommutative probability space $(\mathcal A,\varphi)$ that are zero mean ($\varphi(a_i)=0$), unit variance ($\varphi(a_i^2)=1$), are identically distributed or are subject to some uniform bounds on the higher moments (see hypothesis H3 further on), are `independent' (see hypothesis H4), and pair-wise satisfy the commutation relations
$$a_i a_j\,\,=\,\,s_{i,j}\,a_ja_i$$
where $(s_{i,j})_{i,j\in\mathbb N}$ is some prescribed sequence of elements of  $\{-1,1\}$.
The central question is that of the asymptotic distribution of the sums
$$Z_N:=\frac{1}{\sqrt N}\,\sum_{i=1}^N a_i.\label{eq-speicher}$$

When all of the $(a_i)_{i\in\mathbb N}$ commute, i.e.\ $s(i,j)=1$ for all $i,j\in\mathbb N$, the above setting reduces to the classical case; that is, $Z_N$ converges in distribution (equivalently, in moments) to a standard Gaussian random variable:
$$\lim_{N\to\infty}\varphi(Z_N^k)=\begin{cases}0,&k\text{ odd,}\\(k-1)!!,&k\text{ even,}\end{cases}$$
for all $k\in\mathbb N$.
More generally, given an arbitrary sequence of commutation coefficients $(s(i,j))_{1\leq i<j}$, the moments $\varphi(Z_N^k)$ may not converge as $N\to\infty$. To circumvent the pathological cases, Speicher employed a `stochastic interpolation' step, showing that if the commutation coefficients are drawn i.i.d.\ at random with mean $\E(s(i,j))=q$, almost every sequence of commutation coefficients yields a limit. Furthermore, this limit equals
$$\lim_{N\to\infty}\varphi(Z_N^k)=\begin{cases}0,&k\text{ odd,}\\\sum_{\pi\in\P_2(k)}q^{\cross(\pi)},&k\text{ even,}\end{cases}$$
in which the reader may recognize the moments of the standard $q$-Gaussian distribution of \cite{Bozejko1991}.
As such, Speicher's theorem provides both an independent proof of the positivity of the $q$-commutation relations (see e.g. \cite{Frisch1970,Zagier1992, Bozejko1997} for some related work and historical notes) as well as
a useful method of transferring results from the bosonic/fermionic frameworks to the $q$-Gaussian setting \cite{Biane1997b,Kemp2005}.

Some 20 years following the original result, Blitvi\'c generalized Speicher's theorem by showing that, with some additional care, the commutation `spins' can be extended to real-valued commutation coefficients \cite{BlitvicCLT}.
The corresponding central limits are the $(q,t)$-Gaussian statistics (for $|q|<t$), associated with the commutation relation
  \begin{equation}
a_{q,t}(x)a^\ast_{q,t}(y)-qa^\ast_{q,t}(y)a_{q,t}(x)=\langle x,y\rangle t^N\quad\quad(x,y\in H),\label{com-rel-qt}
\end{equation}
where $a^\ast_{q,t}(x)$ and $a_{q,t}(x)$ now denote the creation and annihilation operators on the $(q,t)$-Fock space of \cite{Blitvic2012}. This two-parameter family also includes the $q$-Gaussian statistics as a special case ($t=1$) and turns out to have connections to a wealth of objects in physics, combinatorics, $q$-series, and other areas (see \cite{Blitvic2012}). The similarity between \eqref{eq-comm} and the above commutation relation \eqref{com-rel-qt} has been observed in \cite{Bozejko15}, and extends down to the moment formulas. Indeed, the moments of the $(q,t)$-Gaussian elements, presently denoted by $\tilde G_{q,t}(x)$, also draw on the combinatorics of crossings and nestings in pair-partitions, as
\begin{align}
&\langle \Omega,\tilde G_{q,t}(x_{2n-1})\ldots \tilde G_{q,t}(x_1)\Omega\rangle_{q,t}=0,\\
&\langle \Omega,\tilde G_{q,t}(x_n)\ldots \tilde G_{q,t}(x_1)\Omega\rangle_{q,t}= \sum_{\pi\in \P_2(2n)}q^{\cross(\pi)}t^{\text{Nest}(\pi)}\prod_{\substack{(i,j)\in\pi}}\langle x_i,x_j\rangle.
 \label{eq-moments-2}
\end{align}
(See Section~\ref{sec-comb} for the relevant combinatorial definitions.) 
The apparent similarities between the noncommutative processes in \cite{Blitvic2012} and \cite{Bozejko15} raise the question of how these may be related. As we now show, both arise from a noncommutative Central Limit Theorem, by generalizing the central argument of \cite{Speicher1992} in two different directions.


\subsection{Main Results}

In the present article, we show that the theorem of Speicher \cite{Speicher1992} can be adapted to recover the Gaussian statistics of type B. Rather than generalizing at the level of the commutation coefficients as in \cite{Blitvic2012}, we retain the commuting/anticommuting structure of \cite{Speicher1992} and instead obtain the desired limits by passing to two (coupled) sequences of elements.

\begin{assump}
Given a noncommutative probability space $(\A,\varphi)$, the main hypotheses for two sequences $(a_i)_{i\in\mathbb N}$ and $(b_i)_{i\in\mathbb N}$ of self-adjoint elements of $\A$ are as follows:
\begin{enumerate}[label=\roman*.]
\item[(H1)] (Vanishing means) For all $i\in\mathbb N$, $\varphi(a_i)=\varphi(b_i)=0$.
\item[(H2)] (Fixed second moments) For all $i\in\mathbb N$, $\varphi(a_i^2)=\varphi(b_i^2)=1$ and $\varphi(a_ib_i)=\rho\in(-1,1)$.

\item[(H3)] (Uniform moment bounds) There exists a sequence $(\gamma_n)_{n\in\mathbb N}$ in $\R_+$ such that for all $n\in\mathbb N$, $i(1),\ldots,i(n)\in\mathbb N$ and  $c_{1,i(1)}\in\{a_{i(1)},b_{i(1)}\},\ldots,c_{n,i(n)}\in\{a_{i(n)},b_{i(n)}\}$, we have
$$\left|\varphi\left(\prod_{j=1}^nc_{j,i(j)}\right)\right|\leq \gamma_n.$$
\item[(H4)] (``Independence'') $\varphi$ factorizes over the naturally ordered products in $\{a_i,b_i\}_{i\in\mathbb N}$. That is, denoting by $\mathcal A_i$ the unital $\ast$-subalgebra generated by $\{a_i,b_i\}$ and letting $g_i\in\mathcal A_i$ ($i\in\mathbb N$), we have
$$\varphi(g_{i(1)}\ldots g_{i(k)})=\varphi(g_{i(1)})\ldots\varphi(g_{i(k)}),$$
for all $k\in\mathbb N$, whenever
$i(1)<i(2)<\ldots<i(k)$.

\item[(H5)] (Commutation relations) There are sequences $(s_{i,j})_{i,j\in\mathbb N}$ and $(r_{i,j})_{i,j\in\mathbb N}$ in $\{-1,1\}$
such that for all $i\neq j$,
$$a_i a_j\,\,=\,\,s_{i,j}\,a_ja_i, \quad\quad b_ib_j\,\,=\,\,s_{i,j}\,b_jb_i,\quad\quad a_ib_j\,\,=\,\,r_{i,j}\,b_ja_i.$$

\end{enumerate}
\label{assump-CLT}
\end{assump}

Observe that the second part of the hypothesis (H2) is partially redundant. Indeed, by the Cauchy-Schwarz inequality,  the unit variances of the elements constrain the range of the parameter $\rho$ to the interval $[-1,1]$, as
\begin{equation*}|\rho|=|\varphi(a_i  b_i)|\leq \sqrt{\varphi(a_i  a_i)\varphi(b_i  b_i)}=1.
\end{equation*}
The reasons for restricting the hypothesis to $\rho\in(-1,1)$ will become apparent shortly.

Compared to \cite{Speicher1992}, there are now two sequences of elements of $\mathcal A$ rather than one. Furthermore, we are interested in the asymptotic distribution of the sums
$$\frac{1}{\sqrt N}\,\sum_{i=1}^N \frac{a_i+b_i}{\sqrt 2}.$$
While $a_i+b_i$ and $a_j+b_j$ may neither commute nor anticommute, with a little care the proof technique of \cite{Speicher1992} remains applicable. We thus show the following.

\begin{theorem}
Let $(\A,\varphi)$ be a noncommutative probability space and
fix $q\in (-1,1)$.
Let
$(a_i)_{i\in\mathbb N}$ and $(b_i)_{i\in\mathbb N}$ be sequences from $\A$ that satisfy the Main Hypotheses
with respect to
commutation coefficients $(s_{i,j})_{1\leq i< j}$ and $(r_{i,j})_{i,j\in\mathbb N}$
drawn i.i.d.\ at random from $\{-1,1\}$, with
\begin{equation}\E(s_{i,j})=\E(r_{i,j})=q.\label{eq-prob}\end{equation}
Set
\begin{equation}S_N:=\frac{1}{\sqrt N}\,\sum_{i=1}^N \frac{a_i+b_i}{\sqrt 2}.\label{eq-SN}\end{equation}
Then, almost surely,
\begin{eqnarray}\lim_{N\to\infty}\varphi(S_N^{2n-1})&=&0,\label{eqLimit1B}\\
\lim_{N\to\infty}\varphi(S_N^{2n})&=&
\sum_{(\pi,f)\in \P_2^B(2n)}\rho^{\text{NB}(\pi,f)}\,q^{\cross(\pi)+2\text{CNB}(\pi,f)}.\label{eqLimit2B}\end{eqnarray} \label{CLT1}

\end{theorem}

Returning to the setting of \cite{Bozejko15}, given a unit vector $e\in H$, a scalar $\alpha\in \R$, and a bounded linear self-adjoint involution $\Pi_0$ on $H$, letting $\rho=\alpha \langle e, \Pi_0 e\rangle$ in \eqref{eqLimit1B} and \eqref{eqLimit2B} recovers the moments of $G_{\alpha,q}(e)$, the Gaussian operator of type B associated with $e$.

As previously observed, the unit variance hypothesis (H2) does not in itself preclude us from considering the boundary cases $|\rho|=1$. While \eqref{eqLimit2B} is not applicable when $|\rho|=1$, as $a_i$ becomes a scalar multiple of $b_i$ and the commutation coefficients $r(i,j)$ and $s(i,j)$ can no longer be drawn independently of one another, the random variables $\varphi(S_N^{2n})$ nevertheless converge to a limit.  The case $\rho=-1$ gives rise to a  degenerate Gaussian element with mean and variance zero, owing to the fact that $a_i+b_i=0$ (hence $S_N=0$), whereas $\rho=1$ recovers the $q$-Gaussian limits with mean zero and variance equal to 2, as $(a_i+b_i)/\sqrt 2=\sqrt 2\, a_i$.

Since $\Pi_0$ is a bounded linear self-adjoint involution on $H$, it follows that $\langle e,\Pi_0  e\rangle\in[-1,1]$ for any unit vector $e\in H$. Since, $\rho$ must also take values in $[-1,1]$, one thus independently recovers the fact that the type B Gaussian elements, $G_{\alpha,q}(x)$ ($x\in H$), are defined for $\alpha\in [-1,1]$, resp. $\alpha\in (-1,1)$ in the strictly positive definite case.


Note that the choice of considering self-adjoint elements with prescribed covariances (H2) is made in order to directly recover the type B Gaussian statistics. To provide asymptotic models for the creation and annihilation operators on the Fock space of type B, one may instead take $\varphi(a_i^\ast a_i)=\varphi(b_i^\ast b_i)=1$, $\varphi(a_i a_i^\ast)=\varphi(b_i b_i^\ast)=\varphi(a_i^2)=\varphi(b_i^2)=\varphi(a_ib_i)=\varphi(a_i^\ast b_i^\ast)=\varphi(a_ib_i^\ast)=0$, and $\varphi(a_i^\ast b_i)=\rho\in(-1,1)$ and consider the mixed moments $\lim_{N\to\infty}\varphi(S_N^{\epsilon_1}\ldots S_N^{\epsilon_k})$ for all $k\in\mathbb N$ and $\epsilon_1,\ldots,\epsilon_k\in\{1,\ast\}$. The given proof adapts easily to give these.

Theorem 1 can be generalized to a type B Gaussian process analogously to \cite{Speicher1992,BlitvicCLT}, with a little additional care. Observe that by \eqref{eq-moments}, depending on the choice of the self-adjoint unitary operator $\Pi_0$, the type B Gaussian elements associated to orthogonal test vectors need not be orthogonal with respect to the vacuum state. For example, fixing an o.n.\ basis $(e_n)$ of $H$ and letting $\Pi_0$ transpose $e_1$ with $e_2$ while leaving the other basis vectors invariant, we have $\varphi_{\alpha,q}(G_{\alpha,q}(e_1)G_{\alpha,q}(e_2))=\alpha$. In order to obtain a noncommutative analogue of Brownian motions in this setting, and, in particular, for $\varphi_{\alpha,q}(G_{\alpha,q}(e_{i}) G_{\alpha,q}(e_{j}))=\langle e_i,(I+\alpha \Pi_0) e_j\rangle$ to vanish whenever $i\neq j$, it suffices for the vectors $(e_n)$ to be orthogonal and the operator $\Pi_0$ to be diagonal with respect to $(e_n)$.
For any $N,k\in\mathbb N$, let 
 \begin{equation}S_{N,k}:=\frac{1}{\sqrt N}\sum_{i=1}^{N}\frac{a_i^{(k)}+b_i^{(k)}}{\sqrt 2},\end{equation}
 where, $(a_i^{(k)})$ and $(b_i^{(k)})$ are pairs of sequences that are uncorrelated in $k$, in the sense that for each $i,j\in\mathbb N$, $$\varphi(a_i^{(k)}a_j^{(k')})=\varphi(b_i^{(k)}a_j^{(k')})=\varphi(b_i^{(k)}b_j^{(k')})=0\quad\text{whenever } k\neq k',$$ and, for each $k\in\mathbb N$, $(a_i^{(k)})$ and $(b_i^{(k)})$ satisfy the Main Hypotheses with covariance $\rho_k=\alpha\langle e_k,\Pi_0e_k \rangle$. We then
recover the moments of the Gaussian process of type B. Namely, the reader may verify that for any choice of $k$ and $i(1),\ldots,i(k)\in\mathbb N$,
\begin{align}&\lim_{N\to\infty}\varphi(S_{N,i(1)} \ldots S_{N,i(k)})=\varphi_{\alpha,q}(G_{\alpha,q}(e_{i(1)}) \ldots G_{\alpha,q}(e_{i(k)})).\end{align}\label{thm-asymptoticB}

Finally, one can construct matrices satisfying the Main Hypotheses with respect to a suitable state, as follows. The result may be referred to as a \emph{Jordan-Wigner transform} (see \cite{Carlen1993,Biane1997b,Kemp2005,BlitvicCLT} for some related constructions), extended to accomodate additional moment hypotheses (H3) and the additional commutativity structure (H5).

\begin{theorem} Fix $|\rho|<1$, as well as the commutation coefficients $\{s_{i,j}\}_{1\leq i<j}$ and $\{r_{i,j}\}_{i,j\in\mathbb N}$. Consider the Hilbert space $\C^2$, vector $(1,0)\in \C^2$, and let $(K,v)=\bigotimes_{i\in\mathbb N} (\C^2,(1,0))$ be the infinite tensor product of $\C^2$ with itself with respect to the constant stabilizing sequence given by $(1,0)$.

For any $x\in\{-1,1\}$, let $\sigma_x,\gamma,\tau\in\mathscr B(\C^2)$ be as follows:
\begin{eqnarray*}&&\sigma_{x}=\left[\begin{array}{cc}1&0\\0& x\end{array}\right],
\quad \gamma=\left[\begin{array}{cc}0&1\\1&0\end{array}\right],\quad \tau=\left[\begin{array}{cc}\rho&\sqrt{1-\rho^2}\\\sqrt{1-\rho^2}&-\rho\end{array}\right].\end{eqnarray*} Furthermore, for $i=1,2\ldots$, consider the following elements of $\B(K)$:
\begin{eqnarray*}&&\zeta_{i}=\sigma_{s(1,i)}\otimes\sigma_{s(2,i)}\otimes\ldots\otimes \sigma_{s(i-1,i)}\otimes \gamma\otimes I\otimes I\otimes\ldots,\\
&&\alpha_{i}=I^{\otimes (i-1)}\otimes \gamma\otimes I\otimes I\otimes \ldots,\\
&&\beta_{i}=\sigma_{s(1,i)r(1,i)}\otimes\ldots\otimes \sigma_{s(i-1,i)r(i-1,i)}\otimes\gamma\otimes\sigma_{s(i+1,i)r(i+1,i)}\otimes\sigma_{s(i+2,i)r(i+2,i)}\otimes\ldots,\\
&&\eta_{i}=I^{\otimes (i-1)}\otimes \tau\otimes I \otimes I\otimes \ldots,\\
&&\theta_{i}=I\otimes I\otimes \ldots,
\end{eqnarray*}
where $I$ denotes the identity on $\C^2$.

Let $\A:=\mathscr B(K^{\otimes 3})$ and let $\varphi$ be the vector state on $\A$ corresponding to the vector $v^{\otimes 3}$. Then, the sequences $(a_i)_{i\in\mathbb N}$ and $(b_i)_{i\in\mathbb N}$ of elements of $\A$ given by
\begin{eqnarray}a_{i}&=&\zeta_{i}\otimes \alpha_{i}\otimes \eta_{i}\\b_{i}&=&\zeta_{i}\otimes \beta_{i}\otimes \theta_{i}\label{eq-JWB}\end{eqnarray}
satisfy the Main Hypotheses with respect to $\varphi$.
\end{theorem}

\section{Combinatorial Objects}\label{sec-comb}

We briefly survey the objects that provide the combinatorial underpinnings of Theorem 1. Let $\P(n)$ denote the collection of partitions of the set $[n]=\{1,\ldots,n\}$. Given $\pi\in\P(n)$, elements of $\pi$ are referred to as the \emph{blocks} of $\pi$. The \emph{size} of a block is the cardinality of the underlying set. (E.g. $\pi=\{\{1,2,4\},\{3\},\{5\}\}\in\P(5)$ is formed by three blocks, one of which has size three and two have size one.)

Two vectors of indices will be declared equivalent if element repetitions occur at same locations in both vectors, namely
for $(i(1),\ldots,i(r)),(j(1),\ldots,j(r))\in [n]^r$,
\begin{align}&(i(1),\ldots,i(r))\sim (j(1),\ldots,j(r))&\iff&\quad \text{ for all } 1\leq k_1<k_2\leq r,&\notag\\
&&&i(k_1)=i(k_2)\text{ iff }j(k_1)=j(k_2).&\label{eq-sim}\end{align}
The equivalence classes of $[n]^{r}$ under ``$\sim$'' are in obvious correspondence with elements of $\P(r)$. (E.g. $(2,2,3,2,4)$ is in the equivalence class corresponding to $\pi=\{(1,2,4),(3),(5)\}$.)

We will be particularly interested in the collection $\P_2(2n)$ of \emph{pair partitions} (aka \emph{pairings}) of $[2n]$, which are partitions whose blocks all have size two. It will be further convenient to represent a pair partition as an ordered list of ordered pairs, that is, $\P_2(2n)\ni\pi=\{(w_1,z_1),\ldots,(w_n,z_n)\}$, where $w_i<z_i$ for $i\in[n]$ and $w_1<\ldots<w_n$. A \emph{pair partition of type B} (in the sense of \cite{Bozejko15}) is a pair $(\pi,f)$, with $\pi\in\P_2(2n)$ (for some $n\in\mathbb N$) and $f:\pi\to\{-1,1\}$ a coloring of the blocks of $\pi$.

The pair partitions of type B will appear with the following combinatorial refinements. For $\pi=\{(w_1,z_1),\ldots,(w_{n},z_{n})\}\in\P_2(2n)$, pairs $(w_i,z_i)$ and $(w_j,z_j)$ are said to \emph{cross} if $w_i<w_j<z_i<z_j$. Let $\Cr(\pi)$ denote the number of pairs of blocks in $\pi$ that cross, namely
\begin{equation}\Cr(\pi):=\#\{(w_i,w_j,z_i,z_j)\mid (w_i,z_i),(w_j,z_j)\in\pi\text{ with } w_i<w_j<z_i<z_j\}.\label{eq-cross}
\end{equation}
Crossings are analogously defined for pair partitions of type B by ignoring the coloring of the blocks, namely $\Cr(\pi,f):=\Cr(\pi)$.

\begin{figure}
\includegraphics[scale=0.5]{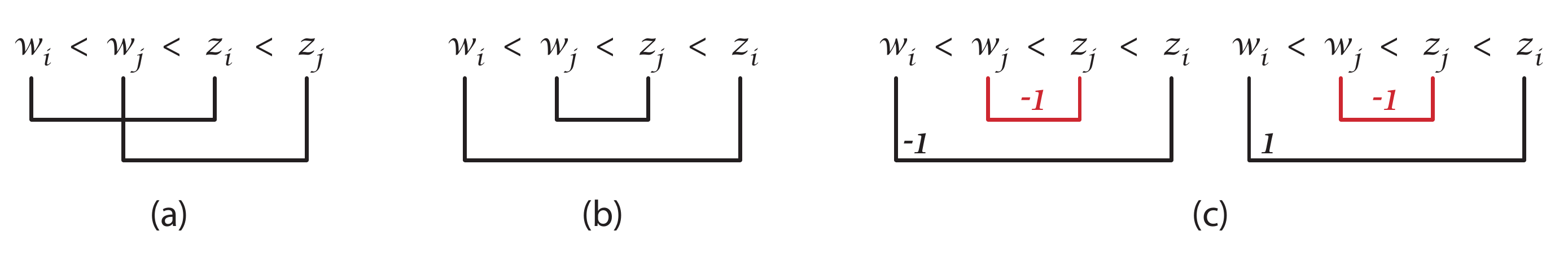}
 \caption{(a) A crossing of two blocks of a pair partition, (b) a nesting of two blocks of a pair partition, (c) an `asymmetric' nesting (CNB) of two blocks of a type B pair partition.}
 \label{fig-comb}
\end{figure}

Analogously, for $\pi=\{(w_1,z_1),\ldots,(w_{n},z_{n})\}\in\P_2(2n)$, pairs $(w_i,z_i)$ and $(w_j,z_j)$ are said to \emph{nest} if $w_i<w_j<z_j<z_i$, and we let
\begin{equation}\Nest(\pi):=\#\{(w_i,w_j,z_i,z_j)\mid (w_i,z_i),(w_j,z_j)\in\pi\text{ with } w_i<w_j<z_j<z_i\}.\label{eq-nest}
\end{equation}
Note that nestings are a natural combinatorial counterpart to crossings; for example, the two combinatorial statistics are equidistributed, in the sense that $$\sum_{\pi\in\P_2(2n)} q^{\Cr(\pi)}=\sum_{\pi\in\P_2(2n)} q^{\Nest(\pi)}.$$ Nestings play a central role in \cite{BlitvicCLT}, arising as a consequence of the passage from commutation signs to real-valued commutation coefficients. In the present case,
when extended to pair partitions of type B, the notion of a nesting ceases to be symmetric, depending instead of the coloring of the blocks. Namely, let
\begin{eqnarray}\CNB(\pi,f):=\#\{(w_i,w_j,z_i,z_j)\mid (w_i,z_i),(w_j,z_j)\in\pi&\text{ with } w_i<w_j<z_j<z_i&\notag\\&\text{ and }f(w_j,z_j)=-1\}.
\end{eqnarray}
(In \cite{Bozejko15}, CNB stands for the number of pairs of ``a covering block'' and a ``negative block''.) Crossings and nestings in pair-partitions and pair partitions of Type B are illustrated in Figure~\ref{fig-comb} .

\section{The Proofs}

Prior to proving the main theorems, several remarks are in order.

\begin{remark}
As $s_{i,j},r_{i,j}\in\{-1,1\}$, the commutation relations are compatible with the $\ast$-structure and the positivity of $\varphi$. In particular, the consistency relations (A), (B) and (C) of \cite{BlitvicCLT} (see p.\ 1464 and p.\ 1469) are automatically met.
\end{remark}
\begin{remark}
As $s_{i,j}^2=r_{i,j}^2=1$, the moment-factorizing hypothesis (H4) above takes on a simpler form when restricted to products whose underlying set partition (see Section~\ref{sec-comb}) is a pair partition. For instance, for the fourth moment $\varphi(a_j  b_j a_i^2)$ with $i<j$, the commutation relations and the independence hypothesis yield
$$\varphi(a_j  b_j a_i^2)=r_{i,j}^2s_{i,j}^2\varphi(a_i^2 a_j  b_j)=\varphi(a_i^2)\varphi(a_j  b_j).$$
More generally, for all $c_i,d_i\in \{a_i,b_i\}$ ($i\in\mathbb N$) and distinct indices $i(1),\ldots,i(k)\in\mathbb N$, we have
$$\varphi(c_{i(1)}d_{i(1)}c_{i(2)}d_{i(2)}\ldots c_{i(k)}d_{i(k)})=\varphi(c_{i(1)}d_{i(1)})\varphi(c_{i(2)}d_{i(2)})\ldots \varphi(c_{i(k)}d_{i(k)}),$$
regardless of the ordering of $i(1),\ldots,i(k)$. In other words, $\varphi$ factorizes over words indexed by pair partitions regardless of the ordering of the pairs. (This is the analogue of Remark 2 of \cite{BlitvicCLT}.) \label{remark-compatibility}
\end{remark}

\begin{proof}[Proof of Theorem 1]
The proof proceeds along analogous lines to \cite{Speicher1992} (see also \cite{BlitvicCLT}). Given any $k\in\mathbb N$, consider the $k$th moment $\varphi(S_N^k)$.
By \eqref{eq-SN}, expressing the product of sums as a sum of products yields the identity
\begin{equation}\varphi(S_N ^k)=\frac{1}{(2N)^{k/2}}\sum_{i(1),\ldots,i(k)\in[N]}\sum_{c}\,\varphi\Big(c(1,i(1)) \ldots c(k,i(k))\Big),\label{eq-SN-prod}\end{equation}
where the inner-most sum is over all
\begin{equation}
 c(1,i(1))\in\left\{a_{i(1)},b_{i(1)}\right\},\ldots,c(k,i(k))\in\left\{a_{i(k)},b_{i(k)}\right\}.\label{eq-c}
\end{equation}
(Note that the heavier indexing notation, compared to \cite{Speicher1992,BlitvicCLT} is due to the fact that we are working with two generators.) Following \cite{Speicher1992}, in order to keep track of which of the elements in the product arise from the same subalgebra, we consider the equivalence class of each $k$-tuple $(i(1),\ldots,i(k))$.
Specifically,
$$\varphi(S_N ^k)=\sum_{\pi\in\P(k)} \frac{1}{(2N)^{k/2}}\sum_{\substack{i(1),\ldots,i(k)\in[N]\text{ s.t. }\\(i(1),\ldots,i(k))\sim\pi}}\sum_{c}\,\varphi\Big(c(1,i(1)) \ldots c(k,i(k))\Big),$$
where we made use of the equivalence relation \eqref{eq-sim} to group together the noncommutative words whose second indices are in the same equivalence class. (For example, $\varphi(a_5\,b_1\,b_5 \, b_7)$ is indexed by the quadruple $(5,1,5,7)$, which is in the equivalence class of the set partition $\pi=\{(1,3),(2),(4)\}$.)

Applying the commutation relations (H5) and the moment-factorizing hypothesis (H4), any expression of the form
\begin{equation}\varphi\Big(c(1,i(1)) \ldots c(k,i(k))\Big)\label{eq-phi}\end{equation}
can now be factorized according to the blocks of the underlying set partition $\pi$. The hypothesis on the vanishing of the means (H1) ensures that no partitions containing a singleton block contribute to \eqref{eq-SN-prod}. (In the previous example, $\varphi(a_5\,b_1\,b_5 \, b_7)=s_{5,1}\,\varphi(a_5\,b_5 )\varphi(b_1)\varphi(b_7)=0$.) Furthermore, a standard counting argument in conjunction with the uniform bounds (H3) ensures that no partitions containing a block of cardinality at least three  contribute to $\lim_{N\to\infty} \varphi(S_N^k)$ (there being too few such partitions compared to the normalizing factor $N^{k/2}$). It follows that
$$\lim_{N\to\infty}\varphi(S_N^{2n-1})=0,$$
while
\begin{align}&\lim_{N\to\infty}\varphi(S_N^{2n})\notag\\
&=\lim_{N\to\infty}\sum_{ \pi\in\P_2(2n)} \frac{1}{2^n\,N^{n}}\sum_{\substack{i(1),\ldots,i(2n)\in[N]\text{ s.t. }\\(i(1),\ldots,i(2n))\sim\pi}}\sum_c\,\varphi\Big(c(1,i(1)) \ldots c(2n,i(2n))\Big),\label{eq-lim-penultimate}\end{align}
where we emphasize that the only contributing set partitions are now the pair partitions.

When $\pi\in\P_2(2n)$, as per Remark~\ref{remark-compatibility}, \eqref{eq-phi} can be expressed as a product of commutation coefficients times a product of (mixed) second moments. Furthermore, it is convenient to notationally distinguish the cases where the second moments arising from the factorization are mixed moments, of the form $\varphi(a_ib_i)$ and $\varphi(b_ia_i)$, as opposed to $\varphi(a_i a_i)$ and $\varphi(b_i b_i)$. This distinction induces a type B pair partition (see Section~\ref{sec-comb}), whose negative blocks are those indexing moments of the former form, while the positive blocks are indexing the moments of the latter. Hence, \eqref{eq-lim-penultimate} can be expressed as
\begin{align}&\lim_{N\to\infty}\varphi(S_N^{2n})\notag\\
&=\lim_{N\to\infty}\sum_{ (\pi,f)\in\P_2^B(2n)} \frac{1}{2^n\,N^{n}}\sum_{\substack{i(1),\ldots,i(2n)\in[N]\text{ s.t. }\\(i(1),\ldots,i(2n))\sim\pi}}\sum_{c\text{ respecting }f}\,\varphi\Big(c(1,i(1))\ldots c(2n,i(2n))\Big),\label{eq-lim-penultimate-2}\end{align}
where the inner-most sum is over all $c(1,i(1))\in\left\{a_{i(1)},b_{i(1)}\right\},\ldots,c(2n,i(2n))\in\left\{a_{i(2n)},b_{i(2n)}\right\}$ such that for every $(w,z)\in\pi$, $c(w,i(w))\neq c(z,i(w))$ if $f(z,w)=-1$ and $c(w,i(w))= c(z,i(w))$ if $f(z,w)=1$. (Recall that $i(w)=i(z)$ for all $(z,w)\in\pi$.)

Since by (H2), $\varphi(a_i^2)=\varphi(b_i^2)=1$ and $\varphi(a_ib_i)=\rho$, it follows that for all indices $i(1),\ldots,i(k)\in\mathbb N$ in the equivalence class of a given $(\pi,f)\in \P_2^B(2n)$,
\begin{equation}\sum_{c\text{ respecting }f}\varphi\Big(c(1,i(1))\ldots c(2n,i(2n))\Big)=\theta_{c(1,i(1))\ldots c(2n,i(2n))}\,\rho^{\text{NB}(\pi,f)}2^{n} ,\end{equation}
where $\theta_{c(1,i(1))\ldots c(2n,i(2n))}$ is the product of the commutation coefficients incurred by commuting the word $c(1,i(1))\ldots c(2n,i(2n))$ into naturally ordered form (see (H4) as well as Remark~\ref{remark-compatibility}), while the factor $2^{n}$ accounts for the fact that there are two choices for each positive block (namely, $a_ia_i$ and $b_ib_i$) and similarly two choices for each negative block (that is, $a_ib_i$ and $b_ia_i$).

To fully recover \eqref{eqLimit2B}, it remains to characterize the family $\theta_{c(1,i(1))\ldots c(2n,i(2n))}$. Since $\pi\in\P_2(2n)$, it suffices to consider how two pairs of elements commute. First, consider moments of the form $\varphi(c_ic_jd_id_j)$, that is, where the indices in the product are in the equivalence class of the pair-partition $\{(1,3),(2,4)\}$ (two pairs that cross). Letting $i\neq j$, by parts (H2), (H4), and (H5) of the Main Hypotheses,
\begin{align*}&\varphi(a_i a_j a_ia_j)  =   s_{i,j}=\varphi(b_i b_j  b_ib_j)\\
&\varphi(a_i a_j b_ib_j)  =  r_{i,j}\,\rho^2=\varphi(b_i b_j a_ia_j)\\
&\varphi(a_i a_j b_ia_j)  = r_{i,j}\,\rho=\varphi(b_i b_j a_ib_j)\\
&\varphi(a_i a_j a_ib_j)  = s_{i,j}\,\rho=\varphi(b_i b_j b_ia_j)\\
&\varphi(a_i b_j a_ia_j)  =  r_{i,j}\,\rho=\varphi(b_i a_j b_ib_j)\\
&\varphi(a_i b_j b_ia_j)  = s_{i,j}\,\rho^2=\varphi(b_i a_j a_ib_j)\\
&\varphi(a_i b_j b_ib_j)  =  s_{i,j}\,\rho=\varphi(b_i a_j a_ia_j)\\
&\varphi(a_i b_j a_ib_j)  =  r_{i,j}=\varphi(b_i a_j b_ia_j).
\end{align*}
Similarly, for the moments of the form $\varphi(c_ic_jd_jd_i)$ where the indices in the product are in the equivalence class of the pair-partition $\{(1,4),(2,3)\}$ (two pairs that nest), we obtain
\begin{align*}&\varphi(a_i  a_j^2a_i)  =   s_{i,j}^2=\varphi(b_i  b_j^2 b_i)\\
&\varphi(a_i a_j  b_jb_i)  =  s_{i,j}r_{i,j}\,\rho^2=\varphi(b_i  b_j  a_ja_i)\\
&\varphi(a_i a_j  b_ja_i)  = s_{i,j}r_{i,j}\,\rho=\varphi(b_i  b_j  a_jb_i)\\
&\varphi(a_i a_j^2 b_i)  = r_{i,j}^2\,\rho=\varphi(b_i  b_j^2 a_i)\\
&\varphi(a_i  b_j  a_ja_i)  =  s_{i,j}r_{i,j}\,\rho=\varphi(b_i  a_j  b_jb_i)\\
&\varphi(a_i  b_j^2a_i)  = r_{i,j}^2=\varphi(b_i  a_j^2 b_i)\\
&\varphi(a_i  b_j^2b_i)  =  s_{i,j}^2\,\rho=\varphi(b_i  a_j^2 a_i)\\
&\varphi(a_i  b_j  a_jb_i)  =  s_{i,j}r_{i,j}\rho^2=\varphi(b_i  a_j  b_ja_i).
\end{align*}
Recalling that the commutation coefficients are drawn at random, consider $\E(\theta_{c(1,i(1))\ldots c(2n,i(2n))})$. By assumption, $\E(r_{i,j})=\E(s_{i,j})=q$ and the distinct commutation coefficients are independent. Hence, each pair of crossing blocks in the underlying pair partition $\pi$, i.e.\ every pair of blocks $(w_\ell,z_\ell),(w_m,z_m)\in\pi$ s.t.\ $w_\ell<w_m<z_\ell<z_m$, contributes a factor of $q$ to $\E(\theta_{c(1,i(1))\ldots c(2n,i(2n))})$, arising as either $\E(s_{i,j})$ or $\E(r_{i,j})$. Furthermore, since $\E(s_{i,j}^2)=\E(r_{i,j}^2)=1$ and $\E(s_{i,j}r_{i,j})=q^2$, only certain nesting pairs contribute nontrivially, namely, those giving rise to terms of the form $s_{i,j}r_{i,j}$. Specifically, each pair of nesting blocks in the underlying pair partition $\pi$, i.e.\ $(w_\ell,z_\ell),(w_m,z_m)\in\pi$ s.t.\ $w_\ell<w_m<z_m<z_\ell$, either contributes a factor of $1$, if the terms of the product indexed by the \emph{inner} block are $(a_{i(w_m)} ,a_{i(w_m)})$ or $(b_{i(w_m)} ,b_{i(w_m)})$, or a factor of $q^2$, if the terms of the product indexed by the \emph{inner} block are $(a_{i(w_m)} ,b_{i(w_m)})$ or $(b_{i(w_m)} ,a_{i(w_m)})$. In other words, it is only those nestings whose inner block (seen as belonging to a partition of type B) is negative which contribute a non-trivial factor, and that factor is $q^2$. Hence,
$$\E(\theta_{c(1,i(1))\ldots c(2n,i(2n))})=q^{\cross(\pi)+2\text{CNB}(\pi,f)},$$ where we recall that $\text{CNB}(\pi,f):=\{(w_\ell,z_\ell),(w_m,z_m)\in\pi\mid w_\ell<w_m<z_m<z_\ell, f(w_m,z_m)=-1\}$ (see Section~\ref{sec-comb}).

We showed that \eqref{eqLimit2B} holds on average. To show that the limit exists and is as stated for a.e.\ sequence of commutation coefficients one may for instance use the Borel-Cantelli lemma. The required estimates can be obtained as in the proof of Lemma 1 of \cite{Speicher1992}, with the calculation following through subject to minor (obvious) modifications.
\end{proof}

\begin{proof}[Proof of Theorem 2]

Start by observing that all the elements are self-adjoint. Let $\varphi_0$ be the vector state on $\mathscr B(\C^2)$ with respect to the vector $(1,0)$ and $\varphi_1$ the vector state on $\mathscr B(K)$ with respect to the vector $v$. In the following, we verify each of the five parts of the Main Hypotheses.

As $\varphi_0(\gamma)=0$, hypothesis (H1)\ is satisfied.

Next, since $\varphi_0(\sigma_x^2)=\varphi_0(\gamma^2)=\varphi_0(\tau^2)=1$, it follows that $\varphi(a_i^2)=\varphi(b_i^2)=1$. Furthermore, since $\varphi_0(\tau)=\rho$, it follows that $\varphi(a_ib_i)=\rho$ and hypothesis (H2)\ is satisfied.

Let $k\in\mathbb N$ and $i_1,\ldots,i_k\in\mathbb N$. Since the commutation coefficients take values in $\{-1,1\}$, it follows that $|\varphi_1(\zeta_{i_1}\ldots\zeta_{i_k})|\leq 1$. Furthermore, letting $\omega_{1,i(1)}\in\{\alpha_{i(1)},\beta_{i(1)}\},\ldots,\omega_{k,i(k)}\in\{\alpha_{i(k)},\beta_{i(k)}\}$, we have $|\varphi_1(\omega_{1,i(1)}\ldots\omega_{k,i(k)})|\leq 1$. Letting instead $\omega_{1,i(1)}\in\{\eta_{i(1)},\theta_{i(1)}\},\ldots,$ $\omega_{k,i(k)}\in\{\eta_{i(k)},\theta_{i(k)}\}$ we similarly obtain $|\varphi_1(\omega_{1,i(1)}\ldots\omega_{k,i(k)})|\leq 1$, as $\tau^2=I$ and $|\rho|<1$. So hypothesis (H3)\ holds.

Now fix $1\leq i_1<i_2<\ldots<i_k$. By linearity of $\varphi$, it suffices to consider the case where each $g_i$ in hypothesis (H4)\ is a product of elements in $\{a_i,b_i\}$. Start by observing that $\sigma_x\sigma_y=\sigma_{xy}$ and $\varphi_0(\sigma_x^m)=1$ for all $m\in\mathbb N$ (hence, $\varphi_0(\sigma_x^m)=\varphi_0(\sigma_x)^m$). Furthermore, $\gamma^2=I$ and $\varphi_0(\sigma_x\gamma^m\sigma_y)=\varphi_0(\sigma_x)\varphi_0(\gamma^m)\varphi_0(\sigma_y)$. (Given the ordering $1\leq i_1<i_2<\ldots<i_k$ and the fact that $\sigma_x\sigma_y=\sigma_{xy}$, these are the only cases to consider.) Hence,  $\varphi_1(\zeta_{i_1}\ldots\zeta_{i_k})=\varphi_1(\zeta_{i_1})\ldots\varphi_1(\zeta_{i_k})$. Similarly, letting $\omega_j\in\{\alpha_j,\beta_j\}$ for all $j\in\{i_1,\ldots i_k\}$, it follows that $\varphi_1(\omega_{i_1}\ldots\omega_{i_k})=\varphi_1(\omega_{i_1})\ldots\varphi_1(\omega_{i_k})$.
Finally, letting instead $\omega_j\in\{\eta_j,\theta_j\}$ for all $j\in\{i_1,\ldots i_k\}$, we observe that $\varphi_0(I^k\tau^m I^n)=\varphi_0(I^k)\varphi_0(\tau^m)\varphi_0(I^n)$. (Given the ordering $1\leq i_1<i_2<\ldots<i_k$, this is the only case to consider.) Hence, $\varphi_1(\omega_{i_1}\ldots\omega_{i_k})=\varphi_1(\omega_{i_1})\ldots\varphi_1(\omega_{i_k})$. Assumption (H4)\ therefore holds.

For the final hypothesis, observe that for $x\in\{-1,1\}$, $\gamma\sigma_{x}=x\sigma_{x}\gamma$. Hence, $\zeta_i\zeta_j=s(i,j)\zeta_j\zeta_i$. Furthermore, $\alpha_i\alpha_j=\alpha_j\alpha_i$ and $\beta_i\beta_j=\beta_j\beta_i$, as $s(i,j)=s(j,i)$ and $s(i,j)^2=1$ for all $i,j\in\mathbb N$. Similarly, $\alpha_i\beta_j=s(i,j)r(i,j)\beta_j\alpha_i$. Finally, $\eta_i\eta_j=\eta_j\eta_i$, $\theta_i\theta_j=\theta_j\theta_i$, and $\eta_i\theta_j=\theta_j\eta_i$. All together, $a_ia_j=s(i,j)a_ja_i$, $b_ib_j=s(i,j)b_jb_i$, and $a_ib_j=r(i,j)b_ja_i$, as required.

For ease of reference, the key ideas behind this construction are summarized in Figure~\ref{fig-matrix}.

\begin{figure}
 \includegraphics[scale=0.32]{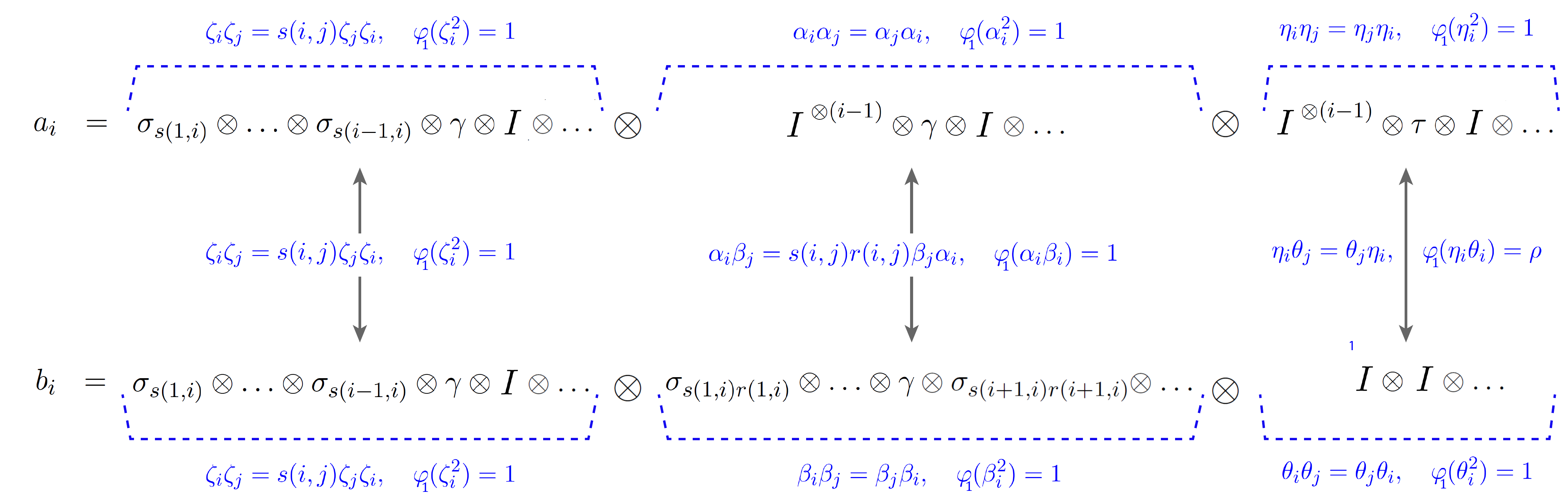}\caption{Matrix models for $a_1,a_2,\ldots, b_1,b_2,\ldots$}\label{fig-matrix}
\end{figure}
\end{proof}

\renewcommand{\abstractname}{Acknowledgement}
\begin{abstract}
\noindent The authors would like to thank Martin Lindsay whose comments greatly improved this manuscript.
\end{abstract}

\bibliographystyle{alpha}
\bibliography{CLT_typeB}

\end{document}